\documentclass[12pt]{article}
\usepackage{amsmath,amsfonts,amsthm,amscd}
\newtheorem{theorem}{Theorem}
\newtheorem{claim}{Claim}
\newtheorem{proposition}{Proposition}
\newtheorem{lemma}{Lemma}
\newtheorem{definition}{Definition}
\newtheorem{conjecture}{Conjecture}
\numberwithin{equation}{section}
\begin{document}
\title{Estimates and Existence Results for some Fully Nonlinear
 Elliptic Equations on Riemannian Manifolds}
\author{Jeff A. Viaclovsky}
\date{February 22, 2001 (revised October 2, 2001)}
\maketitle
\section{Introduction}
We examine the following fully nonlinear partial 
differential equation on a smooth compact $n$-dimensional Riemannian
manifold $(N,g)$
\begin{align}
\label{main}
\sigma_k^{1/k} \left( \nabla^2 u + du \otimes du - \frac{|\nabla u|^2}{2}g
+ S \right) = \psi(x,u) > 0,
%\tag{$1.1_k$}
\end{align}
where $\sigma_k$
%$$\sigma_k(\lambda_1, \dots, \lambda_n) = \sum_{i_1 < \dots < i_k}
%\lambda_{i_1} \cdots \lambda_{i_k},$$
is the $k$th elementary symmetric function of the eigenvalues, 
$S$ is a symmetric tensor, $\nabla$ denotes the
gradient, $\nabla^2$ denotes the Hessian, and 
$du$ is the differential of~$u$. 
%Note that we must use the metric to view the argument 
%of $\sigma_k^{1/k}$ as a symmetric linear
%transformation, and we take the eigenvalues of this map. 
\begin{definition}
Let $(\lambda_1, \dots, \lambda_n) \in \mathbf{R}^n$.
We view the elementary symmetric functions as
functions on $\mathbf{R}^n$
$$\sigma_k(\lambda_1, \dots, \lambda_n) = \sum_{i_1 < \dots < i_k}
\lambda_{i_1} \cdots \lambda_{i_k},$$
and we define
$$\Gamma_k^+ = \mbox{component of } \{\sigma_k > 0\}
\mbox{ containing the positive
cone}.$$
We also define $\Gamma_k^- = - \Gamma_k^+$.
\end{definition}
For a symmetric linear transformation $A : V \rightarrow V$, where
$V$ is an $n$-dimensional inner product space, the
notation $A \in \Gamma_k^{\pm}$ will mean that
the eigenvalues of $A$ lie in the corresponding set.
We note that this notation also makes sense for a symmetric tensor on a
Riemannian manifold. 
%\begin{definition}
%We let $$\Gamma_k^+ =
%\mbox{component of } \{\sigma_k > 0\} \mbox{ containing the
%positive cone}.$$ 
%\end{definition}
%\setcounter{equation}{1}

 We assume the following conditions
\begin{align}
\label{i}
&S \in \Gamma_k^+,
%\tag{$2_k$}
\end{align}
and there exist two constants
$  \underline{\delta} < 0 < \overline{\delta}$ with 
\begin{align}
\label{ii}
\psi(x, \underline{\delta}) < \sigma_k^{1/k}(S) <
\psi(x, \overline{\delta})
\mbox{ for all } x \in N.
%\tag{$3_k$}
\end{align}
For example, we may take $S = g$, and 
$\psi(x,u) = f(x) e^u$, with $f(x) >0$ any smooth positive 
function. We shall see that (\ref{i}) 
is the condition for ellipticity, and 
(\ref{ii}) is the $C^0$ estimate. 
 For equation (\ref{main}) with $1 \leq k \leq n$, we will prove 
\begin{theorem} 
\label{maintheorem}
If $S \in C^{\infty}$, $\psi \in C^{\infty}$,
and both {\em{(\ref{i})}} and {\em{(\ref{ii})}} are satisfied, then there 
exists at least one solution $u \in C^{\infty}(N)$ to {\em{(\ref{main})}}
satisfying $ \underline{\delta} < u < \overline{\delta} $.
\end{theorem}
 In the beautiful paper, \cite{Yanyan1}, Yanyan Li proves
the existence of a solution to
the following equation on a compact Riemannian manifold 
$$ \sigma_k^{1/k}( \nabla^2 u + I) = \psi(x,u) > 0,$$
provided that $N$ has non-negative sectional curvature. 
We would like to emphasize that because of the 
quadratic gradient terms in equation (\ref{main}),
we do not require any curvature assumption in our
existence theorem.

The main part of our proof is the derivation of  
an a priori $C^2$ estimate on \mbox{solutions.} 
The $C^{2,\alpha}$ estimate follows from the 
work of Evans \cite{Evans}, and Krylov \cite{Krylov}
for concave, \mbox{uniformly} elliptic equations.
See also \cite{GT} for an excellent exposition 
of these results.
From these estimates, we obtain the existence theorem 
by applying the degree theory for 
fully nonlinear second order elliptic equations 
developed by Yanyan Li in \cite{Yanyan2}.

  We will also discuss the equations (\ref{main}), when 
$S \in \Gamma_k^-$, the negative cone. By sending 
$u$ to $-u$, we see that the negative case is equivalent 
to the positive cone case of the following equation
\begin{align}
\label{negmain}
\sigma_k^{1/k} \left( \nabla^2 u - du \otimes du + \frac{|\nabla u|^2}{2}g
+ S \right) = \psi(x,u) > 0.
%\tag{$1.4_k$}
\end{align}
In Section \ref{negativesection}, we will show for $\psi(x,u) = 
f(x) e^u$, the $C^1$ estimate still
holds for this equation, but our method for obtaining 
the $C^2$ estimate does not work. We do not know if there 
exists a solution in this case. 
%%%%%%%%%%%%%%%%%%%%%%%%%%%%%%%%%%%%%%%%%%%%%%%%%%%%%%%%%%
\subsection{Conformal Geometry}

   We would also like to point out that (\ref{main})
has geometric origin in conformal geometry; see \cite{Jeff1}. 
Let $(N, g)$ be a Riemannian manifold 
of dimension $n \geq 3$, and we define 
\begin{align*}
A_g \equiv \frac{1}{n-2} \left( Ric - \frac{R}{2(n-1)}g \right),
\end{align*}
where $Ric$ and $R$ are the Ricci tensor and scalar 
curvature of the metric $g$, respectively.
We consider the curvature equation
\begin{align}
\label{curvature equations}
\sigma_k^{1/k} (A_{\tilde{g}}) = \mbox{constant} > 0,
\end{align}
for metrics $ \tilde{g}$ in the conformal class of $g$.
Notice that for $k=1$, the trace, this is 
just the Yamabe equation. 

If we let $\tilde{g} = e^{-2u} g$, then the curvature
equation (\ref{curvature equations})
may be written as the partial differential 
equation (see \cite{Jeff3})
\begin{align}
\label{geometric}
\sigma_k^{1/k} \left( \nabla^2 u + du \otimes du - \frac{|\nabla u|^2}{2}g
+ A_g \right) =  e^{-2u},
\end{align}
where we have normalized the constant to be $1$.
This equation is conformally invariant; see \cite{Jeff2}.

If $A_g \in \Gamma_k^+$, the equation (\ref{geometric})
does not satisfy (\ref{ii}), but our results here 
reduce the compactness question to obtaining 
a $C^0$ estimate on solutions. To this end, for the 
determinant case, we have the following.
Let $ \Omega = \{ \tilde{g} \in [g] : A_{\tilde{g}} \in \Gamma_n^+\}$,
where  $[g]$ denotes the conformal class of $g$, and 
define the conformal invariant
\begin{align}
\sigma([g]) = \underset{\tilde{g} \in \Omega}{ \mbox{inf}}
 \left( \lambda_{max}(A_{\tilde{g}}) D^2 \right),
\end{align}
where $\lambda_{max}(A_{\tilde{g}})$ denotes the maximum 
eigenvalue of the curvature $A_{\tilde{g}}$ on $N$, and $D$ is 
the diameter of $(N, \tilde{g})$.
If $\Omega$ is empty, then define $\sigma([g]) = \infty$. 
\begin{theorem}
\label{detthm}
If $(N,[g]) $ satisifes $\sigma([g]) < \frac{\pi^2}{2}$,  
%\begin{align} 
%\label{geometriccondition}
%\lambda_{max}(A_g) D^2 < 4,
%\end{align}
then there exists $\tilde{g} \in [g]$ satisfying
\begin{align}
\label{zoob}
\mbox{\em{det}}(A_{\tilde{g}}) = 1.
\end{align}
Furthermore, the space of solutions of (\ref{zoob})
is compact.
\end{theorem}
In Section \ref{conformalsection}  will show that, in this 
case, convexity yields a Harnack inequality for solutions which, 
together with a maximum principle argument, produces the 
necessary $C^0$ estimate. To show existence, we use 
a fixed point theorem of Berger (\cite{Berger}), following
an argument from the paper of Delano\"e \cite{Delanoe2}. We will also give 
some examples of manifolds satisfying the 
condition $\sigma([g]) < \frac{\pi^2}{2}$, and
demonstrate that $\sigma( S^n, [g_0]) = \frac{\pi^2}{2}$, 
where $(S^n, g_0)$ is the $n$-sphere with the standard metric. 
Therefore Theorem~\ref{detthm} is analogous to the first step 
in the solution of the Yamabe problem: if the $\sigma$-invariant 
is strictly less than that of the sphere, one has 
existence of solutions and compactness of the space 
of solutions. 

  The case $k=1$, the Yamabe Problem, has been solved 
by Aubin and Schoen (see~\cite{LeeandParker},~\cite{Schoen}), 
and the proof of the $C^0$ estimate for the
Yamabe equation in the locally conformally flat case, 
along with an brief outline of the proof in the general case,
may be found in \cite{Schoen3}. Because of the conformal invariance
of equation (\ref{geometric}), it is reasonable 
to expect that we also have compactness 
for all $k, 1 \leq k \leq n,$ if $(N,g)$ is not conformally 
equivalent to $(S^n, g_0)$:

\begin{conjecture}
\label{conj}
If  $A_g \in \Gamma_k^+$, then there exists 
a conformal deformation $\tilde{g} =e^{2u} g$ such that 
$\sigma_k(A_{\tilde{g}}) = 1$.
Furthermore, if  $(N,g)$ is not conformally equivalent to 
$S^n$ with the standard metric, then the space of solutions 
is compact. 
\end{conjecture}

Again, the results in this paper reduce this compactness 
statement of this conjecture to obtaining $C^0$ estimates on 
solutions. The existence should then follow from a 
suitable topological argument. 
We mention that recently Chang, Gursky and Yang, 
have proved the conjecture for $\sigma_2$ in dimension 
4 (see \cite{CGYnew}).

 Finally, if $A_g \in \Gamma_k^-$, then 
%\begin{align*}
%\sigma_k^{1/k} (C (\tilde{g})) = \mbox{constant} < 0,
%\end{align*}
%for metrics $ \tilde{g}$ in the conformal class of $g$
%(we take the negative $k$th root).
writing $\tilde{g} = e^{2u} g$, 
and normalizing the constant, the equation (\ref{curvature equations})
becomes
\begin{align}
\label{geometric2}
\sigma_k^{1/k} \left( \nabla^2 u - du \otimes du + \frac{|\nabla u|^2}{2}g
- A_g \right) =  e^{2u}.
\end{align}
This is precisely equation (\ref{negmain}) and, 
as mentioned above, from the results in Section~\ref{negativesection},
we have an a priori bound on the $C^1$ norm 
of any solution. We do not know if there 
exists an a priori $C^{2}$ bound for solutions of this equation.

\subsection{Acknowledgements}
  The author is especially grateful to 
Yanyan Li, for numerous helpful 
discussions on fully nonlinear equations, 
and to Pengfei Guan, for pointing out the improvement of
the original argument of Proposition~\ref{aaaa}
which gives the best constant. 
He would also like to thank Alice Chang, 
Phillip Griffiths, Matt Gursky, Karen Uhlenbeck, 
Paul Yang, and Yu Yuan for their interest and 
many useful comments. 

This material is based upon work
supported under a National Science Foundation Postdoctoral
Research Fellowship.
%%%%%%%%%%%%%%%%%%%%%%%%%%%%%%%%%%%%%%%%%%%%%%%%%%%%%%%%%%%%%
\section{Ellipticity}
\label{ellip}
In this section we will show that the equations
(\ref{main}) are elliptic at any solution. 

\begin{definition}
\label{Newtontensor} Let $A : V \rightarrow V$ be 
a symmetric linear transformation
where $V$ is an $n$-dimensional inner product space.
For $0 \leq q \leq n$, the {\em{$q$th Newton transformation}}
associated with $A$ is
\begin{align*}
T_q(A) = \sigma_q(A) \cdot I - \sigma_{q-1}(A)\cdot A + \cdots +
(-1)^qA^q.
\end{align*}
\end{definition}
It is proved in \cite{Reilly} that if $A^i_j$ are the components
of $A$ with respect to some basis of $V$ then
\begin{align}
\label{localNewtontensor} T_q(A)^i_j = \frac{1}{q!} \delta^{i_1
\dots i_q i}_{j_1 \dots j_q j} A_{i_1}^{j_1} \cdots A_{i_q}^{j_q},
\end{align}
where $\delta^{i_1 \dots i_q i}_{j_1 \dots j_q j}$ is the
generalized Kronecker delta symbol, and we are
using the Einstein summation convention. We also have
\begin{align}
\sigma_k(A) = \frac{1}{k!} \delta^{i_1 \dots i_k}_{j_1 \dots j_k}
A_{i_1}^{j_1} \cdots A_{i_k}^{j_k}.
\end{align}
We note that if $A : \mathbf{R} \rightarrow
\mbox{Hom}(V,V)$, then
\begin{align}
\label{dsigmak}
\frac{d}{dt} \sigma_k(A(t)) = T_{k-1}(A(t))^i_j \frac{d}{dt}
A(t)^j_i = T_{k-1} (A(t))^{ij} \frac{d}{dt} A(t)_{ij},
\end{align}
that is, the ($k-1$)-Newton transformation is what we get when we
differentiate $\sigma_k$.

The following proposition describes some important properties
of the sets $\Gamma_k^{+}$.
\begin{proposition}
\label{coneinclusions}
\label{ellipticity}
Each set $\Gamma_k^{+}$ is an open convex cone
with vertex at the origin, and we
have the following sequence of inclusions
$$\Gamma_n^{+} \subset \Gamma_{n-1}^{+} \subset \dots \subset
\Gamma_{1}^{+}.$$
For symmetric linear transformations $A \in \Gamma_k^+$,
$B \in \Gamma_k^+$, and $t \in [0,1]$,
we have the following inequality
\begin{align}
\label{convexity}
\{ \sigma_k( (1-t)A + tB) \}^{1/k}
\geq (1-t) \{ \sigma_k(A) \}^{1/k}
+ t \{ \sigma_k(B) \}^{1/k}.
\end{align}
Furthermore, if $A \in \Gamma_k^+$, then $T_{k-1}(A)$ is
positive definite.
\end{proposition}
The proof of this proposition is standard, and may be found in
\cite{CNSIII} and \cite{Garding}.
Note that by replacing $A$ with $-A$, analogous statements hold 
for $\Gamma_k^-$. Note that the inequality (\ref{convexity})
states that $\sigma_k^{1/k}$ is a {\em{concave}} function 
in $\Gamma_k^+$.

\begin{definition}
\label{admissible}
A function $u \in C^2(N)$ is {\em{positive k-admissible}},
or {\em{negative k-admissible}} if 
\begin{align}
\label{matrix}
\bar{\nabla}^2 u \equiv \nabla^2 u + du \otimes du - 
\frac{|\nabla u|^2}{2}g + S
\end{align}
is everywhere in $\Gamma_k^{+}$ or $\Gamma_k^{-}$, respectively.
\end{definition}
%Note that we are using the metric $g$ to view this
%tensor as an endomorphism
%of the tangent space at any point, and we take the
%eigenvalues of this map. Since the tensor
%is symmetric, these eigenvalues are real.
\begin{proposition}
\label{ellprop}
If $S \in \Gamma_k^+$, then equation {\em{(\ref{main})}} 
is elliptic at any solution.
\end{proposition}
\begin{proof}
Since $N$ is compact, at a minimum of the solution $u$
we have
\begin{align*}
%\label{minpoint}
\sigma_k^{1/k} \left( \nabla^2 u (p)+ S(p)
 \right) = \psi(p,u(p)) > 0,
\end{align*} 
with $\nabla^2u$ positive semidefinite. 
From Proposition \ref{ellipticity}, we then have,
at the minimum point, $\bar{\nabla}^2 u $ is in $\Gamma_k^{+}$.
Therefore since the cones are connected,
by continuity we have $u$ is positive $k$-admissible.
A similar argument holds in the negative 
$k$-admissible case. 
  
\begin{claim}
If we make the conformal change
of metric $\tilde{g} = e^{-2u} g$, then for any function h,  
\begin{align}
\label{hessian}
{\nabla}^2_{\tilde{g}} (h)= \nabla^2_{g}(h) + du \otimes dh + dh \otimes du
- \langle du, dh \rangle g.  
\end{align}
where $\nabla^2_g(h)$ is the Hessian of $h$ with respect
to the metric $g$, and $\nabla^2_{\tilde{g}}(h)$ is taken 
with respect to $\tilde{g}$. 
\end{claim}
\begin{proof}
We have for the Christoffel symbols (see \cite{Besse})
\begin{align*}
\tilde{\Gamma}^l_{ij} = \Gamma^l_{ij} - u_i \delta^l_j
-u_j \delta^l_i + g_{ij}g^{lr}u_r.
\end{align*}
Therefore
\begin{align*}
({\nabla}^2_{\tilde{g}} h)_{ij} &= h_{ij} - \tilde{\Gamma}^l_{ij}h_l\\
&= h_{ij} - ( \Gamma^l_{ij} - u_i \delta^l_j
-u_j \delta^l_i + g_{ij}g^{lr}u_r )h_l\\ 
& = (\nabla^2_g h)_{ij} + u_i h_j + u_j h_i - g^{lr}u_r h_l g_{ij}.
\end{align*}
\end{proof}

We let $$F [ u, \nabla u, \nabla^2 u] =
\sigma_k^{1/k} \left( \nabla^2 u +du \otimes du 
- \frac{|\nabla u |^2}{2} g_0 + S \right) - \psi(x,u).$$
From (\ref{dsigmak}) and  (\ref{hessian}), 
we see that the linearization at the solution $u$ in the 
direction $h$ is given by
\begin{align}
\label{lin}
 F'[u, \nabla u, \nabla^2 u](h) =
 \sigma_k( \bar{\nabla}^2 u )^{\frac{1-k}{k}}
 T_{k-1} ( \bar{\nabla}^2 u )^{ij}
(\nabla^2_{\tilde{g}} h)_{ij} 
- \psi_u h.
\end{align}
Since $\bar{\nabla}^2 u$
is in $\Gamma_k^+$, from Proposition
\ref{ellipticity}, we are done. 
\end{proof}
%%%%%%%%%%%%%%%%%%%%%%%%%%%%%%%%%%%%%%%%%%%%%%%%%%%%%%%%%%%%%
\section{$C^0$ estimate}
\label{C0}
In this section, we present the necessary $C^0$ estimate
which will be required in the existence proof. 
We will give the general argument, and then also 
an easier argument in the case that  
$\psi(x,u) = f(x) e^u$.
%$\psi(x,s) : M \times \mathbf{R} \rightarrow 
%\mathbf{R}^+$ is a proper function.
In order to apply the maximum principle, we need
to rewrite the equation as follows. We let $w = e^u$,
and the equations (\ref{main}) become
\begin{align}
\label{expmain}
\sigma_k^{1/k} \left( \frac{1}{w}  
\nabla^2 w - \frac{1}{w^2} \frac{|\nabla w|^2}{2}g
+ S \right) = \psi(x, \mbox{ln} w),
\end{align}
\begin{proposition}
\label{slp}
If $w_0$ is positive $k$-admissible, and $w_1$ is
positive $k$-admissible, then $(1-t)w_0 + t w_1$
is positive $k$-admissible for $t \in [0,1]$. 
\end{proposition}
\begin{proof}
By positive $k$-admissible, we mean that the matrix 
$$\bar{\nabla}^2 w \equiv
w \nabla^2 w - \frac{|\nabla w|^2}{2}g + w^2 S$$
is in $\Gamma_k^+$. The multiple of $w$ is irrelevant, 
since $w = e^u > 0$.  
Letting $w_t(x) = (1-t) w_0(x) + t w_1(x)$,
we must show that $\bar{\nabla}^2 w_t \in \Gamma_k^+$, 
i.e., $F_k$ is elliptic at $w_t$ for 
$t \in [0,1]$.
We have
\newcommand{\go}{|\nabla w_0|^2}
\newcommand{\gl}{|\nabla w_1|^2}
\begin{align*}
\bar{\nabla}^2 w_t &= w_t \nabla^2 w_t - \frac{|\nabla w_t|^2}{2} g
+ ( (1-t) w_0 + t w_1)^2 S\\
&= ((1-t)w_0 + tw_1)((1-t)\nabla^2 w_0 + t \nabla^2 w_1) - \frac{|\nabla
  ((1-t)w_0+ t w_1)|^2}{2} g\\
& \ \ \ \ \ \ \ \ \ \ \ \ + (1-t)^2w_0^2 S + 2 t(1-t)w_0 w_1 S + t^2 w_1^2 S\\
& = (1-t)^2 w_0 \nabla^2 w_0 + t^2 w_1 \nabla^2 w_1 + t(1-t)(w_0 
\nabla^2 w_1 + w_1  \nabla^2 w_0)\\
& \ \ \ \ \ \ \ \ - \left( (1-t)^2\frac{\go}{2} + t(1-t)\nabla w_1
 \cdot \nabla w_0 + t^2 \frac{\gl}{2} \right) g\\
& \ \ \ \ \ \ \ \ \ \ \ \ + (1-t)^2w_0^2 S + 2 t(1-t)w_0 w_1 S + t^2 w_1^2 S\\
& = (1-t)^2 \bar{\nabla}^2 w_0 + t^2  \bar{\nabla}^2 w_1 + t(1-t)
\biggl( \frac{w_0}{w_1} \Bigl( w_1 \nabla^2 w_1 - \frac{\gl}{2} g +
\frac{\gl}{2} g \Bigr) \\
& + \frac{w_1}{w_0} \Bigl( w_0 \nabla^2 w_0 - \frac{\go}{2} 
g + \frac{\go}{2} g \Bigr)
- (\nabla w_1 \cdot \nabla w_0) g \biggr)+ 2 t(1-t)w_0 w_1 S \\
& = (1-t) \Bigl( (1-t) \bar{\nabla}^2 w_0 + t \frac{w_1}{w_0} 
\bar{\nabla}^2 w_0 \Bigr)
+ t \Bigl( t \bar{\nabla}^2 w_1 + (1-t) \frac{w_0}{w_1} 
\bar{\nabla}^2 w_1 \Bigr)\\
& \ \ \ \ \ \ \ + \frac{t(1-t)}{2w_0 w_1} \Bigl( w_0^2 \gl +
 w_1^2 \go - 2 w_0 \nabla w_1
\cdot w_1 \nabla w_0 \Bigr) g \\
& = (1-t) \Bigl( (1-t) \bar{\nabla}^2 w_0 + 
t \frac{w_1}{w_0} \bar{\nabla}^2 w_0 \Bigr)
+ t \Bigl( t \bar{\nabla}^2 w_1 + (1-t) \frac{w_1}{w_0} 
\bar{\nabla}^2 w_1 \Bigr)\\
& + \frac{t(1-t)}{2w_0 w_1} \Bigl( |w_0 \nabla w_1 - w_1
\nabla w_0|^2 \Bigr) g .
\end{align*}
From Proposition \ref{coneinclusions}, the first two terms together 
are in $\Gamma_k^+$. The last term is a non-negative multiple
of the identity, so again using Proposition 
\ref{coneinclusions}, we are done. 
\end{proof}
\begin{proposition}
\label{C0estimate}
Suppose $S \in C^0$, $\psi \in C^1$, and both 
{\em{(\ref{i})}} and {\em{(\ref{ii})}} are satisfied.
Then any $C^2$ solution $u$ of {\em{(\ref{main})}} with 
$ \underline{\delta} \leq u \leq \overline{\delta}$
satifies $ \underline{\delta} < u < \overline{\delta}.$
\end{proposition}
\begin{proof}
 Assume we have a solution $u$ of (\ref{main}),
 with $\underline{\delta} \leq u$. 
We let 
\begin{align*}
F[w] = \sigma_k^{1/k} \left( \frac{1}{w}  
\nabla^2 w - \frac{1}{w^2} \frac{|\nabla w|^2}{2}g
+ S \right) - \psi(x, \mbox{ln} w).
\end{align*}
Then letting $w = e^u$, the function $w - e^{\underline{\delta}}
 \geq 0$ satisfies
\begin{align*}
L (w - e^{\underline{\delta}} ) = F[w] - F[e^{\underline{\delta}}] 
= 0 - \sigma_k^{1/k}(S) + \psi(x, \underline{\delta}) < 0,
\end{align*}
where $L$ is a linear elliptic operator (this follows 
from Proposition \ref{slp}, see \cite{GT}, Chapter~17), so by the maximum 
principle, we have $e^{\underline{\delta}} < w$, that is, 
$\underline{\delta} < u$.
The proof of the strict upper inequality is similar. 
\end{proof}
\noindent
{\bf{Remark.}} Why did we change to $w = e^u$ in the above argument? 
A computation similar to that of the proof 
of Proposition \ref{slp} shows that the original equation
(\ref{main})
is elliptic along the straight line path only if
$k \leq n/2$. We are just using a different straight line
path in order to apply the maximum principle. 

  In the case that 
%$\psi(x,s) : M \times \mathbf{R} \rightarrow 
%\mathbf{R}^+$ is a proper function, 
$\psi(x,u) = f(x) e^u$,
we present an 
alternative, more elementary derivation of the 
$C^0$ estimate. 
\begin{lemma}
Let $A$ and $B$ be symmetric $n \times n$ matrices. 
Assume that $A$ is positive semi-definite, $B \in \Gamma_k^+$, 
and $A + B \in \Gamma_k^+$. Then 
$$ \sigma_k(A + B) \geq \sigma_k (B).$$
If $A$ is negative semi-definite, then 
$$ \sigma_k(A + B) \leq \sigma_k (B).$$
\end{lemma}
\begin{proof}
Let $F(t) = \sigma_k (tA + B) - \sigma_k(B)$ for $t \in [0,1]$.
Note that from convexity of the cone $\Gamma_k^+$, 
we have $ t(A+ B) + (1-t) B = tA + B \in \Gamma_k^+$. 
Using (\ref{dsigmak}), we have 
\begin{align*}
F'(t) = T_{k-1}(tA + B)^{ij} A_{ij}  \geq 0,
\end{align*}
since $T_{k-1}(tA + B)$ is positive definite from 
Proposition \ref{coneinclusions}.
Therefore $F(t)$ is non-decreasing, and $F(0) = 0$, so we have
$F(1) = \sigma_k(A + B) - \sigma_k(B) \geq 0.$
The negative case is similar. 
\end{proof}
\begin{proposition} 
\label{easyC0}
Suppose $S \in C^0$ satisfies {\em{(\ref{i})}}. If 
%$\psi(x,s) : M \times \mathbf{R} \rightarrow 
%\mathbf{R}^+$ is a proper function 
$\psi(x,u) = f(x) e^u$, for $f(x) > 0$ a positive
$C^0$ function, then there 
exist constants $  \underline{\delta} < 0 < \overline{\delta}$
depending only upon $f$, $S$ and $k$, 
such that for any solution $u(x)$ of {\em{(\ref{main})}}, we have 
$\underline{\delta} < u(x) < \overline{\delta}$.
\end{proposition}
\begin{proof}
Since $N$ is compact, at a minimum of the function $u(x)$
we have
\begin{align*}
%\label{minpoint}
\sigma_k^{1/k} \left( \nabla^2 u (p)+ S(p)
 \right) = f(p) e^{u(p)}
\end{align*} 
with $\nabla^2u(p)$ positive semidefinite. 
From the lemma we have 
\begin{align*}
\sigma_k^{1/k} (S(p)) \leq f(p) e^{u(p)},
\end{align*}
and certainly we can choose $\underline{\delta}$ such that
\begin{align*}
u(x) \geq u(p) \geq \mbox{ln} \left( \frac{ \sigma_k^{1/k}(S(p))}{f(p)} \right)
\geq \mbox{ln} \left( \underset{x \in N}{ \mbox{min}} 
\frac{ \sigma_k^{1/k}(S(x))}{f(x)} \right) > \underline{\delta}.
\end{align*}
Similarly, if the maximum of $u(x)$ is at $q \in N$, we
can choose  $\overline{\delta}$ such that
\begin{align*}
u(x) \leq u(q) \leq \mbox{ln} \left( \frac{ \sigma_k^{1/k}(S(q))}{f(q)} \right)
\leq \mbox{ln} \left( \underset{x \in N}{ \mbox{max}} 
\frac{ \sigma_k^{1/k}(S(x))}{f(x)} \right) < \overline{\delta} .
\end{align*}
\end{proof}
%%%%%%%%%%%%%%%%%%%%%%%%%%%%%%%%%%%%%%%%%%%%%%%%%%%%%%%%%%%%%%%%%5
\section{$C^1$ estimate}
\label{C1}
\begin{proposition}
\label{C1estimate}
Suppose $S \in C^1$, $\psi \in C^1$, 
{\em{(\ref{i})}} is satisfied, and  
$u$ is a $C^3$ solution of {\em{(\ref{main})}}
satisfying $ \underline{\delta} \leq u(x) \leq \overline{\delta}$.
Then there exists a constant $C_1$ depending only upon 
$S, \psi, \underline{\delta}, \overline{\delta},$
and $k$
such that
\begin{align*}
|\nabla u|_{C^0} \leq C_1.
\end{align*}
\end{proposition}
We consider the following function
\begin{align*}
h = \left( 1 + \frac{|\nabla u|^2}{2} \right) e^{\phi(u)},
\end{align*}
where $\phi : \mathbf{R} \rightarrow \mathbf{R}$ is a 
function of the form 
\begin{align*}
\phi(s) = c_1(c_2 - s)^p.
\end{align*}
The constants $c_1, c_2,$ and $p$ will be chosen later. 
We will estimate the maximum value of the 
function $h$, and this will give us the gradient 
estimate. 

  Since $N$ is compact, and $h$ is 
continuous, we suppose the maximum of $h$ occurs and a point 
$p \in N$. We take a normal coordinate system 
$(x^1, \dots, x^n)$ at $p$. Then we have 
$g_{ij}(p) = \delta_{ij}$, and $\Gamma^i_{jk}(p) = 0$, 
where $g = g_{ij} dx^i dx^j$, and $\Gamma^i_{jk}$
is the Christoffel symbol (see \cite{Besse}). 

 Locally, we may write $h$ as 
\begin{align*}
h = \left( 1 + \frac{1}{2} g^{lm}u_{l}u_{m} \right)
e^{\phi(u)} = v e^{\phi(u)}.
\end{align*}
In a neighborhood of $p$, differentiating $h$ in the $x^i$
direction we have
\begin{align}
\notag
\partial_i h &= h_i = \frac{1}{2} \partial_i ( g^{lm}u_l u_m)
e^{\phi(u)} + v e^{\phi(u)}\phi'(u)u_i\\
\label{goog}
& = \frac{1}{2} \partial_i ( g^{lm})u_l u_m e^{\phi(u)} 
 + g^{lm} \partial_i(u_l) u_m e^{\phi(u)} + v e^{\phi(u)}\phi'(u)u_i 
\end{align}
Since in a normal coordinate system, the first 
derivatives of the metric vanish at $p$, 
and since $p$ is a maximum for $h$, 
evaluating (\ref{goog}) at $p$, we have
\begin{align}
\label{hip}
u_{li}u_l = - v \phi'(u) u_i.
\end{align}
Next we differentiate (\ref{goog}) in the $x^j$ direction. 
Since $p$ is a maximum, $\partial_j \partial_i h = h_{ij}$ is negative 
semidefinite, and we get (at $p$)
\begin{align*}
0 \gg h_{ij} &= \frac{1}{2} \partial_j \partial_i g^{lm} u_l u_m
e^{\phi(u)} + u_{lij}u_l e^{\phi(u)} + u_{li}u_{lj} e^{\phi(u)}
+ u_{li}u_{l} e^{\phi(u)} \phi'(u) u_j \\
&+ v_j e^{\phi(u)} \phi'(u) u_i
 + v e^{\phi(u)} (\phi'(u))^2 u_i u_j + v e^{\phi(u)} \phi''(u)u_j u_i
+ v e^{\phi(u)} \phi'(u) u_{ij}
\end{align*} 
Next we note that $v_j = u_{lj} u_l$, and 
using (\ref{hip}), we have 
\begin{align*}
0 \gg h_{ij} &= \frac{1}{2} \partial_j \partial_i g^{lm} u_l u_m
e^{\phi(u)} + u_{lij}u_l e^{\phi(u)} + u_{li}u_{lj} e^{\phi(u)}\\
&+ ( \phi''(u) - \phi'(u)^2) v e^{\phi(u)} u_i u_j 
 + v e^{\phi(u)} \phi'(u) u_{ij}
\end{align*} 
Next we divide by $v e^{\phi(u)}$, sum 
with $T_{k-1}(\bar{\nabla}^2 u)^{ij}$ (which is
positive definite and symmetric), and we have
the inequality
\begin{align}
\label{inequality1}
0 \geq \frac{1}{2v} T_{k-1}^{ij} \partial_i \partial_j
g^{lm}u_l u_m + \frac{1}{v} T_{k-1}^{ij} u_{lij}u_l
+ ( \phi''(u) - \phi'(u)^2)T_{k-1}^{ij} u_i u_j
+ \phi'(u) T_{k-1}^{ij} u_{ij},
\end{align}
since $u_{li} u_{lj}$ is positive semidefinite, 
and we abbreviate $T_{k-1}^{ij} = T_{k-1}(\bar{\nabla}^2 u)^{ij}$,
where $ \bar{\nabla}^2 u$ is the notation in (\ref{matrix}) above.

 We will use equation (\ref{main}) to replace the $u_{ij}$ term
with lower order terms,  and then differentiate equation
(\ref{main}) in order to replace the $u_{lij}$ term with 
lower order terms. Writing equation (\ref{main}) with 
respect to our local coordinate system, we have
\begin{align}
\label{equation}
\sigma_k^{1/k} \biggl( g^{lj} \Bigl( u_{ij} - u_r \Gamma^r_{ij} + u_i u_j
- \frac{1}{2} (g^{r_1 r_2} u_{r_1} u_{r_2} )g_{ij} + S_{ij}
 \Bigl) \biggl) = \psi (x,u).
\end{align}
Note that the $g^{lj}$ term is present since we need 
to raise an index on the tensor before we apply $\sigma_k^{1/k}$.

  For a symmetric matrix $A$ , we have the formula (see \cite{Reilly})
\begin{align*}
T_{k-1}(A)^{ij} A_{ij} = k \sigma_k(A).
\end{align*}
Using this, and equation (\ref{equation}), we have at $p$,
\begin{align}
\notag
T_{k-1}^{ij}u_{ij}& = T_{k-1}^{ij} \left( u_{ij} +u_i u_j 
- \frac{| \nabla u|^2}{2} \delta_{ij} + S_{ij}
-u_i u_j +\frac{| \nabla u|^2}{2} \delta_{ij} - S_{ij} \right)\\
\notag &= k \sigma_k + T_{k-1}^{ij} \left( -u_i u_j +\frac{| \nabla u|^2}{2} 
\delta_{ij} - S_{ij} \right)\\
\label{zero}
& = k \psi(x,u)^k + T_{k-1}^{ij} \left( -u_i u_j +\frac{| \nabla u|^2}{2} 
\delta_{ij} - S_{ij}\right) 
\end{align}
Next we take $m$ with $1 \leq m \leq n$, and apply $\partial_m$
to (\ref{equation})
\begin{align}
\begin{split}
\label{fullonederivative}
& \sigma_k^{\frac{1-k}{k}} T_{k-1}^{il} \bigg(
\partial_m g^{lj} (\bar{\nabla}^2 u)_{ij}
 + g^{lj} \Big( u_{ijm} - u_{rm} \Gamma^r_{ij}
- u_r \partial_m \Gamma^r_{ij} + u_{im} u_j + u_i u_{mj}\\
&- \frac{1}{2}(\partial_m g^{r_1 r_2}) u_{r_1} u_{r_2} g_{ij}
- g^{r_1 r_2} u_{r_1m} u_{r_2} g_{ij} - 
\frac{1}{2} g^{r_1 r_2} u_{r_1} u_{r_2} \partial_m g_{ij} 
+ \partial_m S_{ij} \Big) \bigg)\\
& \hspace{10mm} = \frac{ \partial \psi}{\partial x^m} + 
\frac{ \partial \psi}{\partial u } u_m.
\end{split}
\end{align}
We evaluate the above expression at $p$, and we obtain
\begin{align}
\label{onederivative}
\psi^{1-k} T_{k-1}^{ij} \Big(
 u_{ijm} - u_r \partial_m \Gamma^r_{ij} + 2 u_{im} u_j 
 - u_{rm} u_{r} \delta_{ij} + \partial_m S_{ij} \Big)
= \psi_m + \psi_u u_m.
\end{align}
We then sum with $u_m$, and using (\ref{hip})
we have the following formula
%\begin{align}
%T_{k-1}^{ij} u_{ijm}u_m
%= T_{k-1}^{ij} \Big(
% u_m u_r \partial_m \Gamma^r_{ij} - 2 u_m u_{im} u_j 
% + u_m u_{rm} u_{r} \delta_{ij} + \partial_m S_{ij} \Big)
%- \psi^{k-1}(u_m \psi_m + \psi_u |\nabla u|^2).
%\end{align}
%Now we use (\ref{hip}) and we arrive at
\begin{align}
\label{first}
\notag T_{k-1}^{ij} u_{ijm}u_m
= T_{k-1}^{ij} \Big(
 u_m u_r \partial_m \Gamma^r_{ij} + 2 v \phi'(u) u_i u_j 
 - v \phi'(u) |\nabla u|^2 \delta_{ij} + u_m \partial_m S_{ij} \Big)\\
+ \psi^{k-1}(u_m \psi_m + \psi_u |\nabla u|^2).
\end{align}

Substituting (\ref{zero}) and (\ref{first}) into 
(\ref{inequality1}), we arrive at the inequality
\begin{align}
\notag
0  \geq \frac{1}{2v} &T_{k-1}^{ij} \partial_i \partial_j
g^{lm}u_l u_m \\
\notag
&+ \frac{1}{v}
 T_{k-1}^{ij} \Big(
 u_l u_r \partial_l \Gamma^r_{ij} + 2 v \phi'(u) u_i u_j 
 - v \phi'(u) |\nabla u|^2 \delta_{ij} + u_l \partial_l S_{ij} \Big)\\
\notag
& + \frac{\psi^{k-1}}{v} \Big( u_m \psi_m + \psi_u |\nabla u|^2 \Big)
+ ( \phi''(u) - \phi'(u)^2)T_{k-1}^{ij} u_i u_j\\
\label{mess}
& + k \phi'(u) \psi(x,u)^k + \phi'(u) T_{k-1}^{ij} \left( -u_i u_j
 +\frac{| \nabla u|^2}{2} \delta_{ij} - S_{ij}\right)  
\end{align}
\begin{lemma}
At $p$, in normal coordinates, we have
\begin{align*}
 \sum_{l,m} (\partial_i \partial_j
g^{lm} + 2 \partial_l \Gamma^m_{ij}) u_l u_m 
= 2 \sum_{l,m} R_{iljm}u_l u_m,
\end{align*}
where $R_{iljm}$ are the components of the Riemann
curvature tensor of $g$ {\em{(see~\cite{Besse})}}.
\end{lemma}
\begin{proof}
The metric is parallel, so we have
\begin{align*}
0 = \nabla_j g^{lm} = \partial_j g^{lm} + \Gamma^l_{jr} g^{rm}
+ \Gamma^m_{jr} g^{lr}.
\end{align*}
Therefore we have, at $p$,
\begin{align*}
0= \partial_i \partial_j g^{lm} + \partial_i \Gamma^l_{jr} \delta^{rm}
+ \partial_i \Gamma^m_{jr} \delta^{lr} = 
 \partial_i \partial_j g^{lm} + \partial_i \Gamma^l_{jm} 
+ \partial_i \Gamma^m_{jl}. 
\end{align*}
Using this, we have
\begin{align*}
 \sum_{l,m} (\partial_i \partial_j
g^{lm} + 2 \partial_l \Gamma^m_{ij}) u_l u_m & =
 \sum_{l,m} ( - \partial_i \Gamma^l_{jm} 
- \partial_i \Gamma^m_{jl}+ 2 \partial_l \Gamma^m_{ij}) u_l u_m \\
& =  2 \sum_{l,m} ( - \partial_i \Gamma^m_{lj} 
+ \partial_l \Gamma^m_{ij}) u_l u_m \\
& =  2 \sum_{l,m} R_{iljm}  u_l u_m.
\end{align*}
\end{proof}
Using the lemma, and collecting terms in (\ref{mess}), we arrive at
\begin{align}
\notag
- \psi^{k-1} \Big( \frac{u_m}{v} \psi_m & + 
\psi_u \frac{|\nabla u|^2}{v} \Big) 
 - k \phi'(u) \psi(x,u)^k
 \geq 
%\frac{1}{2v} T_{k-1}^{ij} \Big( \partial_i \partial_j
%g^{lm} + 2 \partial_l \Gamma^m_{ij} \Big) u_l u_m   \\
 \Big( \phi''(u) - \phi'(u)^2 + \phi'(u) \Big) T_{k-1}^{ij} u_i u_j\\
\label{goodone}
& + T_{k-1}^{ij} \left( 
 R_{iljm} \frac{u_l u_m}{v}- \phi'(u) \frac{| \nabla u|^2}{2} 
\delta_{ij} - \phi'(u) S_{ij}
+ \frac{u_l}{v} \partial_l S_{ij} \right).  
\end{align}
 Now we will choose $\phi(s)$. 
\begin{lemma}
\label{choosephi}
Assume that $ \underline{\delta} < s < \overline{\delta}$.
Then we may choose constants $c_1, c_2,$ and $p$
depending only upon $\underline{\delta}$, and $\overline{\delta}$.
so that $\phi(s) = c_1(c_2 -s)^p$ satisfies
\begin{align}
\label{phi'}
 \phi'(s) < 0,
\end{align}
and 
\begin{align}
\label{phi''}
 \phi''(s) - \phi'(s)^2 + \phi'(s) > 0.
\end{align}
\end{lemma}
\begin{proof}
We have 
\begin{align*}
\phi'(s) = -pc_1(c_2 -s)^{p-1},
\end{align*}
and
\begin{align*}
\phi''(s) = p(p-1) c_1 (c_2 -s)^{p-2}.
\end{align*}
To satisfy (\ref{phi'}) we need $c_1 > 0$, $p>0$, and
$c_2 > s$. So choose $c_2 > \overline{\delta}$.
Next we have
\begin{align*}
\phi''(s) - \phi'(s)^2 + \phi'(s) &= 
 p(p-1) c_1 (c_2 -s)^{p-2} - ( pc_1(c_2 -s)^{p-1})^2
- p c_1(c_2 -s)^{p-1}\\
&= pc_1 (c_2 - s)^{p-2} \Big( (p-1) - pc_1 (c_2 -s)^p
- (c_2 -s) \Big).
\end{align*}
Now choose 
\begin{align*}
c_1 = \frac{1}{p^2 \cdot \mbox{ max} \{ (c_2 -s)^p \} },
\end{align*}
and $p$ so large that 
\begin{align*}
\overline{\delta} < c_2 < \underline{\delta} + p -1 - \frac{1}{p}.
\end{align*}
Then we have
\begin{align*}
\phi''(s) - \phi'(s)^2 + \phi'(s) &\geq 
 \frac{1}{ p \cdot \mbox{ max} \{ (c_2 -s)^p \} } (c_2 -s)^{p-2}
\Big( p - 1 - \frac{1}{p} - c_2 + s \Big)\\
& > \frac{1}{ p \cdot \mbox{ max} \{ (c_2 -s)^p \} } (c_2 -s)^{p-2}
( -\underline{\delta} + s ) >0.
\end{align*} 
\end{proof}
With $\phi(s)$ chosen as above, we let 
$$\epsilon_1 = - \mbox{max}\{ \phi'(s) \},$$ and
$$\epsilon_2 = \mbox{min} \{ \phi''(s) - \phi'(s)^2 + \phi'(s) \}.$$  
From the inequality (\ref{goodone}), we have
\begin{align}
\label{thingy}
C \geq 
%\frac{1}{2v} T_{k-1}^{ij} \Big( \partial_i \partial_j
%g^{lm} + 2 \partial_l \Gamma^m_{ij} \Big) u_l u_m   \\
 \epsilon_2  T_{k-1}^{ij} u_i u_j
 + T_{k-1}^{ij} \left( 
 R_{iljm} \frac{u_l u_m}{v} + \epsilon_1 \frac{| \nabla u|^2}{2} 
\delta_{ij} + \phi'(u) S_{ij}
+ \frac{u_l}{v} \partial_l S_{ij} \right),  
\end{align}
where in this equation, and in what follows,
$C$ is a constant depending on 
$\overline{\delta}, \underline{\delta},$ and $\psi$.

Without loss of generality, assume that $\bar{\nabla}^2 u$
is diagonal at $p$. 
Now if for some $i$, a diagonal entry of the
matrix in parenthesis above satisfies
\begin{align*}
 R_{ilim} \frac{u_l u_m}{v} + \epsilon_1 \frac{| \nabla u|^2}{2} 
 + \phi'(u) S_{ii}
+ \frac{u_l}{v} \partial_l S_{ii} < 1,
\end{align*}
then we have the gradient bound. So we may assume that 
\begin{align*}
 R_{ilim} \frac{u_l u_m}{v} + \epsilon_1 \frac{| \nabla u|^2}{2} 
 + \phi'(u) S_{ii}
+ \frac{u_l}{v} \partial_l S_{ii} \geq 1,
\end{align*}
for all $i$.
From the inequality (\ref{thingy}), we conclude that
\begin{align}
\label{thingy2}
C \geq 
%\frac{1}{2v} T_{k-1}^{ij} \Big( \partial_i \partial_j
%g^{lm} + 2 \partial_l \Gamma^m_{ij} \Big) u_l u_m   \\
 \epsilon_2   \sum_{i} T_{k-1}^{ii} u_i^2
 + \sum_{i} T_{k-1}^{ii}.
\end{align}
Noting that 
\begin{align}
\label{Newtontrace}
\sum_{i} T_{k-1}^{ii} = (n-k+1) \sigma_{k-1},
\end{align} 
(see \cite{Reilly}) we deduce that
\begin{align*}
  \sigma_{k-1} \leq C.
\end{align*}
\begin{proposition}
\label{compactness}
Let $k \geq 2$, and $A \in \Gamma_k^+$ be a symmetric 
linear transformation. If $0 < c_1 \leq \sigma_k(A)$,
and $\sigma_{k-1}(A) \leq c_2$, then we have a bound on 
the eigenvalues of $A$, that is, $|\lambda(A)| \leq C$,
where $C$ depends only on $c_1$ and $c_2$.
\end{proposition}
\begin{proof}
 The proof may be found in \cite{Yanyan1}.
\end{proof}
Using this result, if $k \geq 2$, we see that 
\begin{align*}
|\lambda| \leq C,
\end{align*}
and since $T_{k-1}$ is positive definite, this
implies
\begin{align*}
T_{k-1}^{ii} \geq \frac{1}{C}>0, \mbox{ for } i = 1 \dots n.  
\end{align*}
Equation (\ref{thingy2}) then implies that
\begin{align*}
|\nabla u|^2 \leq C.
\end{align*}
Note that in the case $k=1$, we do not require the proposition since
$T_0^{ij} = \delta^{ij}$, and therefore 
(\ref{thingy2}) gives the gradient bound.
%%%%%%%%%%%%%%%%%%%%%%%%%%%%%%%%%%%%%%%%%%%%%%%%%%%%%%%%%%%%%%%%%
\section{$C^2$ estimate}
\label{C2}
\begin{proposition}
Suppose $S \in C^2$, $\psi \in C^2$, 
{\em{(\ref{i})}} is satisfied,   
$u$ is a $C^4$ solution of {\em{(\ref{main})}}
satisfying $ \underline{\delta} \leq u(x) \leq \overline{\delta}$,
and $|\nabla u| < C_1$.
Then there exists a constant $C_2$ depending only upon 
$S, \psi, \underline{\delta}, \overline{\delta}, C_1$,
and $k$
such that
\begin{align*}
|\nabla^2 u|_{C^0} \leq C_2.
\end{align*}
\end{proposition}
Let $S(TN)$ denote the unit tangent bundle of $N$, and 
we consider the following function $w: S(TN) \mapsto \mathbf{R}$,
\begin{align*}
w(e_p) = ( \nabla^2 u + du \otimes du + S)(e_p, e_p).
\end{align*}
Since $S(TN)$ is compact, 
let $w$ have a maximum at the vector $\tilde{e}_p$.
We use normal coordinates at $p$, and by rotating,
assume that the tensor is diagonal at $p$, 
and without loss of generality, we may assume 
that $\tilde{e}_p = {\partial}/{\partial x^1}$.
  
 We let $\tilde{w}$ denote the function defined 
in a neighborhood of $p$
\begin{align*}
\tilde{w}(x) &= ( \nabla^2 u + du \otimes du + S)( {\partial}/{\partial
 x^1}, {\partial}/{\partial x^1} )\\
& = (\nabla^2 u)_{11} + u_1^2 + S_{11}\\
& = u_{11} - \Gamma^l_{11}u_l + u_1^2 + S_{11}.
\end{align*}
Differentiating in the $i$th coordinate direction, we obtain
\begin{align}
\label{wix}
\tilde{w}_i = u_{11i} - \partial_i \Gamma^l_{11} u_l
- \Gamma^l_{11} u_{li} + 2 u_1 u_{1i} + \partial_i S_{11}.
\end{align}
The function $\tilde{w}(x)$ has a maximum at $p$,
so evaluating (\ref{wix}) at $p$, we obtain
\begin{align}
\label{wip}
u_{11i} =  \partial_i \Gamma^l_{11} u_l
- 2 u_1 u_{1i} - \partial_i S_{11}.
\end{align}
Next we differentiate (\ref{wix}) in the $x^j$ direction. 
Since $p$ is a maximum, $\partial_j \partial_i \tilde{w} = \tilde{w}_{ij}$
is negative semidefinite, and we get (at $p$)
\begin{align*}
0 \gg \tilde{w}_{ij} = u_{11ij} - \partial_i \partial_j
\Gamma^l_{11} u_l - \partial_i \Gamma^l_{11} u_{lj}
- \partial_j \Gamma^l_{11} u_{li}
+ 2u_{1j} u_{1i} + 2u_1 u_{1ij} + \partial_i \partial_j
S_{11}.
\end{align*}
We sum with $T_{k-1}(\bar{\nabla}^2 u)^{ij}$
(which is positive definite and symmetric), and we have
the inequality
\begin{align}
\label{testineq}
\begin{split}
0 \geq T_{k-1}^{ij} u_{11ij} -T_{k-1}^{ij} \partial_i \partial_j
\Gamma^l_{11} u_l - 2T_{k-1}^{ij} \partial_i \Gamma^l_{11} u_{lj}\\
+ 2T_{k-1}^{ij} u_{1j} u_{1i} + 2T_{k-1}^{ij} u_1 u_{1ij}
+ T_{k-1}^{ij} \partial_i \partial_j S_{11}.
\end{split}
\end{align}
We will use (\ref{onederivative}) to replace the fifth term, and 
we will differentiate equation (\ref{main}) twice 
to replace the first term. 

We recall that the equation is 
$$\sigma_k^{1/k}( \bar{\nabla}^2 u) = \psi(x,u).$$
To simplify notation, write $ f = \sigma_k^{1/k}$.
Differentiating once in the $x^1$ direction, we had
(equation (\ref{fullonederivative}))
\begin{align*}
\frac{\partial f}{\partial r_{ij}}( \partial_1 (\bar{\nabla}^2 u)^i_j)
= \psi_1 + \psi_u u_1.
\end{align*}
Differentiating twice, we obtain
\begin{align*}
&\partial_1 \left(\frac{\partial f}{\partial r_{ij}} \right)
( \partial_1 (\bar{\nabla}^2 u)^i_j)
+ \frac{\partial f}{\partial r_{ij}}
(\partial_1 \partial_1 (\bar{\nabla}^2 u)^i_j)\\
&= \Big( \frac{\partial^2 f}{\partial r_{ij} \partial r_{lm}} \Big)
( \partial_1 (\bar{\nabla}^2 u)^l_m)
( \partial_1 (\bar{\nabla}^2 u)^i_j)
+ \frac{\partial f}{\partial r_{ij}}
(\partial_1 \partial_1 (\bar{\nabla}^2 u)^i_j)\\
&= \psi_{11} + 2 \psi_{1u}u_1 + \psi_{uu}u_1^2 + \psi_u u_{11}.
\end{align*}
Since $\sigma_k^{1/k}$ is concave in $\Gamma_k^+$,
we have the inequality
\begin{align}
\label{ineq1}
T_{k-1}^{ij} \Big( \partial_1 \partial_1 (\bar{\nabla}^2 u)^i_j \Big)
\geq \psi^{k-1}( \psi_{11} + 2 \psi_{1u}u_1 + \psi_{uu}u_1^2 + 
\psi_u u_{11}).
\end{align}
From formula (\ref{fullonederivative}), we can
expand the left hand side, and evaluate at $p$ to get
\begin{align*}
T_{k-1}^{ij} \Big( \partial_1 \partial_1 (\bar{\nabla}^2 u)^i_j \Big)
&= T_{k-1}^{ij} \Big( \partial_1 \partial_1 g^{jl} (\bar{\nabla}^2 u)_{il}\Big)
+ T_{k-1}^{ij} \Big( u_{ij11} - 2u_{r1} \partial_1 \Gamma^r_{ij}
- u_r \partial_1 \partial_1 \Gamma^r_{ij}\\
&+ 2 u_j u_{i11} + 2u_{i1} u_{j1} - \frac{1}{2}  \partial_1 \partial_1 
g^{r_1 r_2} u_{r_1} u_{r_2} \delta_{ij}
- (u_{r11} u_r + u_{r1}u_{r1}) \delta_{ij}\\
&- \frac{1}{2} |\nabla u|^2   \partial_1 \partial_1 (g_{ij})
+  \partial_1 \partial_1 S_{ij} \Big).
\end{align*}
From (\ref{wip}) we can replace terms of the form 
$u_{11i}$ and we have
\begin{align}
\label{bigmess}
\begin{split}
T_{k-1}^{ij} \Big( \partial_1 \partial_1 (\bar{\nabla}^2 u)^i_j \Big)
&= T_{k-1}^{ij} \Big( \partial_1 \partial_1 g^{jl} (\bar{\nabla}^2 u)_{il}\Big)
+ T_{k-1}^{ij} \Big( u_{ij11} - 2u_{r1} \partial_1 \Gamma^r_{ij}
- u_r \partial_1 \partial_1 \Gamma^r_{ij}\\
&+ 2 u_j ( \partial_i \Gamma^l_{11} u_l
- 2 u_1 u_{1i} - \partial_i S_{11})+ 2u_{i1} u_{j1} - \frac{1}{2}  \partial_1 \partial_1 g^{r_1 r_2} u_{r_1} u_{r_2} \delta_{ij}\\
&- (u_r(\partial_r \Gamma^l_{11} u_l
- 2 u_1 u_{1r} - \partial_r S_{11}) + u_{r1}u_{r1}) \delta_{ij}\\
& \ \ \ \ \ \ \ - \frac{1}{2} |\nabla u|^2   \partial_1 \partial_1 (g_{ij})
+  \partial_1 \partial_1 S_{ij} \Big).
\end{split}
\end{align}
Substituting (\ref{bigmess}) in (\ref{ineq1}), we have
\begin{align}
\label{bigmess2}
\begin{split}
T_{k-1}^{ij} u_{ij11} \geq - & T_{k-1}^{ij} ( \partial_1 
\partial_1 g^{jl} (\bar{\nabla}^2 u)_{il})\\
& + T_{k-1}^{ij} ( 2u_{r1} \partial_1 \Gamma^r_{ij}
+ u_r \partial_1 \partial_1 \Gamma^r_{ij} - 2u_{i1} u_{j1}) \\
& +2T_{k-1}^{ij} ( -u_j \partial_i \Gamma^l_{11} u_l
+ 2 u_j u_1 u_{1i} + u_j \partial_i S_{11})\\
& + T_{k-1}^{ij} ( u_r \partial_r \Gamma^l_{11} u_l
- 2 u_1 u_r u_{1r} - u_r \partial_r S_{11} + u_{r1}u_{r1})
\delta_{ij}\\
&+ T_{k-1}^{ij} \Big( \frac{1}{2}  \partial_1 \partial_1 g^{r_1 r_2} 
u_{r_1} u_{r_2} \delta_{ij} + \frac{1}{2} |\nabla u|^2 
 \partial_1 \partial_1 g_{ij}
-  \partial_1 \partial_1 S_{ij} \Big)\\
& +  \psi^{k-1}( \psi_{11} + 2 \psi_{1u}u_1 + \psi_{uu}u_1^2 + 
\psi_u u_{11}).
\end{split}
\end{align}
Next we will substitute inequality (\ref{bigmess2})
into (\ref{testineq}). Note that the fourth term on the 
right hand side of (\ref{bigmess2}) will cancel the fourth 
term in (\ref{testineq}). We also use equation 
(\ref{onederivative}) to replace the fifth term 
in (\ref{testineq}). We have 
\begin{align*}
0 \geq - & T_{k-1}^{ij} \Big( \partial_1 
\partial_1 g^{jl} (\bar{\nabla}^2 u)_{il} \Big) 
+ T_{k-1}^{ij} \Big( 2u_{r1} \partial_1 \Gamma^r_{ij}
+ u_r \partial_1 \partial_1 \Gamma^r_{ij} \Big)\\
& +2T_{k-1}^{ij} \Big( -u_j \partial_i \Gamma^l_{11} u_l
+ \boxed{2 u_j u_1 u_{1i}} + u_j \partial_i S_{11} \Big)\\
& + T_{k-1}^{ij} \Big( u_r \partial_r \Gamma^l_{11} u_l
- \boxed{2 u_1 u_r u_{1r}} - u_r \partial_r S_{11} + u_{r1}u_{r1} \Big)
\delta_{ij}\\
&+ T_{k-1}^{ij} \Big( \frac{1}{2}  \partial_1 \partial_1 g^{r_1 r_2} 
u_{r_1} u_{r_2} \delta_{ij} + \frac{1}{2} |\nabla u|^2 
 \partial_1 \partial_1 g_{ij}
-  \partial_1 \partial_1 S_{ij} \Big)\\
& +  \psi^{k-1} \Big( \psi_{11} + 2 \psi_{1u}u_1 + \psi_{uu}u_1^2 + 
\psi_u u_{11} \Big)
 -T_{k-1}^{ij} \partial_i \partial_j
\Gamma^l_{11} u_l - 2T_{k-1}^{ij} \partial_i \Gamma^l_{11} u_{lj}\\
& + 2  T_{k-1}^{ij} \Big( u_1 u_r \partial_1 \Gamma^r_{ij} - 
\boxed{2 u_1 u_{i1} u_j} 
 + \boxed{u_1 u_{r1} u_{r} \delta_{ij}} - u_1 \partial_1 S_{ij} \Big)\\
& + 2\psi^{k-1} \Big( \psi_1 u_1 + \psi_u u_1^2 \Big)
+ T_{k-1}^{ij} \partial_i \partial_j S_{11}
\end{align*}
 Note that the boxed terms cancel. Using the bounds on 
lower order quantities, the above simplifies to 
\begin{align}
\label{second}
\begin{split}
C + C \sum_i T_{k-1}^{ii} & \geq  -  T_{k-1}^{ij} \Big( \partial_1 
\partial_1 g^{jl} (\bar{\nabla}^2 u)_{il} \Big) 
+ 2T_{k-1}^{ij} \Big( u_{r1} \partial_1 \Gamma^r_{ij} \Big)\\
& + T_{k-1}^{ij} u_{r1} u_{r1} \delta_{ij}
+   \psi^{k-1}\psi_u u_{11} 
- 2T_{k-1}^{ij} \partial_i \Gamma^l_{11} u_{lj}.
\end{split}
\end{align}
In this equation, and in what follows, $C$ is a
constant depending only on 
$S, \psi, \underline{\delta}, \overline{\delta}, C_1$,
and $k$.

  The next step is to rewrite the second derivative 
terms in terms of $ \bar{\nabla}^2 u$. 
To further simplify notation, we let $\bar{u}_{ij} = (\bar{\nabla}^2 u)_{ij}$.
We have
\begin{align*}
 u_{ij} = \bar{u}_{ij} - u_i u_j + (|\nabla u|^2/2) \delta_{ij}
- S_{ij}.
\end{align*}
Substituting this into (\ref{second}), we obtain
\begin{align*}
C + C \sum_i T_{k-1}^{ii} & \geq  -  T_{k-1}^{ij} \Big( \partial_1 
\partial_1 g^{jl} \bar{u}_{il} \Big)\\
& + 2T_{k-1}^{ij} \Big(  \bar{u}_{r1} - u_r u_1 + (|\nabla u|^2/2) 
\delta_{r1} - S_{r1} \Big) \partial_1 \Gamma^r_{ij}\\
& + \Big( \sum_i T_{k-1}^{ii} \Big)
  \sum_r \Big(  \bar{u}_{r1} - u_r u_1 + (|\nabla u|^2/2) 
\delta_{r1} - S_{r1} \Big)^2\\
&+   \psi^{k-1}\psi_u \Big(  \bar{u}_{11} - u_1^2 + (|\nabla u|^2/2) 
 - S_{11} \Big)\\
&- 2T_{k-1}^{ij} \partial_i \Gamma^l_{11} 
\Big( \bar{u}_{lj} - u_l u_j + (|\nabla u|^2/2) \delta_{lj}
- S_{lj} \Big).
\end{align*}
Next we use the fact that $\bar{u}_{ij}$ is diagonal, 
and absorbing lower order terms we obtain
\begin{align*}
 C + C \sum_i T_{k-1}^{ii} + C \bar{u}_{11} \sum_i T_{k-1}^{ii} \geq  
 -  &\sum_i T_{k-1}^{ii} ( \partial_1 \partial_1 g^{ii} \bar{u}_{ii})
- 2 \sum_i T_{k-1}^{ii} \partial_i \Gamma^i_{11} \bar{u}_{ii}\\
+ & \bar{u}_{11}^2 \sum_i T_{k-1}^{ii} + C \bar{u}_{11}.
\end{align*}
We estimate the first two terms on the right hand side
\begin{align*}
 \sum_i T_{k-1}^{ii} ( \partial_1 \partial_1 g^{ii} \bar{u}_{ii})
+ 2 \sum_i T_{k-1}^{ii} \partial_i \Gamma^i_{11} \bar{u}_{ii}
= \sum_i T_{k-1}^{ii} ( R_{1i1i} \bar{u}_{ii}) \leq
C \underset{i}{\mbox{ max}} | \bar{u}_{ii}| \sum_i T_{k-1}^{ii}.
\end{align*}
Since we are in the cone $\Gamma_k^+$, the trace is
positive by Proposition \ref{coneinclusions},
and since $\bar{u}_{11}$ is the largest eigenvalue,
we have
\begin{align*}
 |\bar{u}_{ii}| \leq (n-1) \bar{u}_{11}, \ \ \ \ i = 1 \dots n.
\end{align*}
Therefore we obtain
\begin{align}
\label{coolineq}
 C + C \bar{u}_{11} + C \sum_i T_{k-1}^{ii} 
+ C \bar{u}_{11} \sum_i T_{k-1}^{ii} \geq  
\bar{u}_{11}^2 \sum_i T_{k-1}^{ii}.
\end{align}
Dividing by $\bar{u}_{11}^2$ and using (\ref{Newtontrace}), we obtain
\begin{align}
\label{secder}
\begin{split}
\sigma_{k-1}  \leq 
\Big( \frac{C}{ \bar{u}_{11}^2} + \frac{C}{\bar{u}_{11}}
\Big) \sigma_{k-1} + \frac{C}{ \bar{u}_{11}^2} 
+ \frac{C}{\bar{u}_{11}}. 
\end{split}
\end{align}
If  
\begin{align*}
\frac{C}{ \bar{u}_{11}^2} + \frac{C}{\bar{u}_{11}}  \geq \frac{1}{2},
\end{align*}
then we have the necessary eigenvalue bound.
So we may assume that
\begin{align*}
\frac{C}{ \bar{u}_{11}^2} + \frac{C}{\bar{u}_{11}}
  \leq \frac{1}{2},
\end{align*}
and substitution into inequality (\ref{secder}) yields
\begin{align*}
\frac{1}{2} \sigma_{k-1}  \leq 
\frac{C}{ \bar{u}_{11}^2} 
+ \frac{C}{\bar{u}_{11}}. 
\end{align*}
Without loss of generality we may assume that $\bar{u}_{11} \geq 1$,
and from the above inequality we obtain
\begin{align*}
\sigma_{k-1} \leq C,
\end{align*}
which by Proposition \ref{compactness} yields the 
eigenvalue bound in the case $k \geq 2$. 
In the case $k=1$, (\ref{secder}) already gives the 
eigenvalue estimate. 
%%%%%%%%%%%%%%%%%%%%%%%%%%%%%%%%%%%%%%%%%%%%%%%%%%%%%%%%%%%%%%
\section{Existence}
\label{existence}
We now prove Theorem \ref{maintheorem}.
The main tool will be the degree theory for fully 
nonlinear second order elliptic equations as
developed in \cite{Yanyan2}.
We consider for $t \in [0,1]$ the family of equations
\begin{align}
\label{maint}
t \sigma_k^{1/k} + (1-t) \sigma_1 = t \psi(x,u) + (1-t)
\sigma_1(S) e^u,
\end{align}
where we abbreviate $\sigma_k^{1/k} = \sigma_k^{1/k}( 
\bar{\nabla}^2u)$.
Note that at $t=0$, the equation is 
\begin{align*}
 \Delta u + \frac{2-n}{2} |\nabla u|^2 + \sigma_1(S) = \sigma_1(S)
e^u.
\end{align*}
From the maximum principle, $u = 0$ is the unique solution. 
\begin{proposition}
For any $t \in [0,1]$, any $C^2$ solution $u^t$ of {\em{(\ref{maint})}} 
with $ \underline{\delta} \leq u \leq \overline{\delta}$
satifies $ \underline{\delta} < u < \overline{\delta}.$
\end{proposition}
\begin{proof}
From assumption (\ref{ii}), we have 
\begin{align*}
t \psi(x, \underline{\delta}) + (1-t) \sigma_1(S) e^{\underline{\delta}}
 < t \sigma_k^{1/k}(S) + (1-t)\sigma_1(S) 
< t\psi(x, \overline{\delta}) + (1-t) \sigma_1(S) e^{\overline{\delta}},
\end{align*}
therefore the proof of Proposition \ref{C0estimate} applies.
\end{proof}
\begin{proposition}
Let $ t \in [0,1]$, and $u^t$ be a solution to {\em{(\ref{maint})}}
with $ \underline{\delta} < u^t < \overline{\delta}$.
Then 
$$ \parallel \! u^t  \! \! \parallel_{C^2} < C,$$
for some constant $C$ independent of $t$.
\end{proposition}
\begin{proof}
 We let $f_t = t \sigma_k^{1/k} + (1-t) \sigma_1.$ 
Define 
$$ \Gamma_{k,t}^+ \equiv   
 \mbox{component of } \{ t \sigma_{k}^{1/k} + (1-t) \sigma_1 > 0\}
\mbox{ containing the positive cone}.$$
Then all of the estimates in the previous sections
hold with $\sigma_k^{1/k}$ replaced by $f_t$, 
and $\Gamma_k^+$ replaced by $\Gamma_{k,t}^+$,
and it is then not difficult to see that we can 
choose $C$ independent of $t$, since the $C^0$
estimate holds uniformly.

\end{proof}
  The above estimate yields uniform ellipticity, and since  
our equation is convex with respect to the second derivative
variables,  by the work of 
Evans \cite{Evans}, and Krylov \cite{Krylov} mentioned
in the introduction, and standard elliptic theory,
there exists a constant $M$ independent of $t$ such that
$$ \parallel \! u^t  \! \! \parallel_{C^{4,\alpha}} < M.$$
Define the subset $\mathcal{O}_t$ of $C^{4,\alpha}$ by
\begin{align*}
\mathcal{O}_t \equiv  & \{ \underline{\delta} < u^t < \overline{\delta} \}
\cap \{ \parallel \! u^t  \! \! \parallel_{C^{4,\alpha}} < M \}\\
& \cap \{ \bar{\nabla}^2 u^t \in \Gamma_{k,t}^+ \}
\cap \{  t \sigma_{k}^{1/k} + (1-t) \sigma_1
> t \delta_0 + (1-t) \sigma_1(S) e^{\underline{\delta}} \},
\end{align*}
where $\delta_0$ is a constant chosen such that 
$\psi(x,s) > \delta_0$ for 
$ \underline{\delta} < s < \overline{\delta}$. 
Define $F_t : C^{4, \alpha} \rightarrow C^{2,\alpha}$ by
\begin{align*}
F_t(u) = t \sigma_k^{1/k} (\bar{\nabla}^2 u)
+ (1-t) \sigma_1 (\bar{\nabla}^2 u)
- t \psi(x,u) - (1-t) \sigma_1(S) e^u.
\end{align*}
There are no solutions of the equation $F_t(u) =0$ on 
$\partial \mathcal{O}_t$,
so the degree of $F_t$ is well-defined and independent
of $t$. 
As mentioned above, there is a unique 
solution at $t=0$. Furthermore, the linearization 
at $u=0$ is invertible. Therefore 
\begin{align*}
\mbox{deg}(F_0, \mathcal{O}_0, 0) = \pm 1,
\end{align*}
and since the degree is independent of $t$, we have
\begin{align*}
\mbox{deg}(F_1, \mathcal{O}_1, 0) = \pm 1,
\end{align*}
and we conclude that (\ref{main}) has a solution
in $\mathcal{O}_1$.

 Note that in the case $\psi(x,u) = f(x)e^u$, we can 
avoid using degree theory since the linearization
is invertible, and the existence follows by using 
the continuity method.
%%%%%%%%%%%%%%%%%%%%%%%%%%%%%%%%%%%%%%%%%%%%%%%%%%%%%%%%%%%%%
\section{The negative cone equation}
\label{negativesection}
 As mentioned in the introduction, the negative cone 
case of (\ref{main}) is equivalent to the 
positive cone case of equation (\ref{negmain}). 
We no longer necessarily have ellipticity along the straight 
line path for this equation (the proof of Proposition 
\ref{slp} does not work for this equation), so we just consider 
the equation
\begin{align}
\label{specific}
\sigma_k^{1/k} \left( \nabla^2 u - du \otimes du + \frac{|\nabla u|^2}{2}g
+ S \right) = f(x) e^u> 0.
%\tag{$7.1_k$}
\end{align}
In this section we will show that we still have the 
$C^0$ and $C^1$ estimate for solutions of this equation. 
The proof Proposition \ref{easyC0} still 
works for this equation, so we have
\begin{proposition} 
Suppose $S \in C^0$ satisfies {\em{(\ref{i})}}. Then there 
exist constants $  \underline{\delta} < 0 < \overline{\delta}$
depending only upon $f$ and $S$, 
such that for any solution $u(x)$ of {\em{(\ref{specific})}}, we have 
$\underline{\delta} < u(x) < \overline{\delta}$.
\end{proposition}
The $C^1$ estimate also holds, with appropriate 
modifications to the proof of Proposition \ref{C1estimate}.
\begin{proposition}
Suppose $S \in C^1$, {\em{(\ref{i})}} is satisfied, and  
$u$ is a $C^3$ solution of {\em{(\ref{specific})}}
satisfying $ \underline{\delta} \leq u(x) \leq \overline{\delta}$.
Then there exists a constant $C_1$ depending only upon 
$S, \psi, \underline{\delta}, \overline{\delta},$
and $k$ such that
\begin{align*}
|\nabla u|_{C^0} \leq C_1.
\end{align*}
\end{proposition}
\begin{proof}
We consider the following function
\begin{align*}
h = \left( 1 + \frac{|\nabla u|^2}{2} \right) e^{\phi(u)},
\end{align*}
where ${\phi} : \mathbf{R} \rightarrow \mathbf{R}$ is a 
function of the form 
\begin{align*}
{\phi}(s) = c_1(c_2 + s)^p.
\end{align*}
The proof procedes exactly as before, but we end 
up with the following analogue of equation
(\ref{goodone}) 
\begin{align}
\notag
 -\psi^{k-1} \Big( \frac{u_m}{v} \psi_m & + 
\psi_u \frac{|\nabla u|^2}{v} \Big) 
 - k \phi'(u) \psi(x,u)^k
 \geq 
 \Big( \phi''(u) - \phi'(u)^2 - \phi'(u) \Big) T_{k-1}^{ij} u_i u_j\\
\label{neggoodone}
& + T_{k-1}^{ij} \left( 
 R_{iljm} \frac{u_l u_m}{v} + \phi'(u) \frac{| \nabla u|^2}{2} 
\delta_{ij} - \phi'(u) S_{ij}
+ \frac{u_l}{v} \partial_l S_{ij} \right).  
\end{align}
\begin{lemma}
Assume that $ \underline{\delta} < s < \overline{\delta}$.
Then we may choose constants $c_1, c_2,$ and $p$
depending only upon $\underline{\delta}$, and $\overline{\delta}$.
so that ${\phi}(s) = c_1(c_2 + s)^p$ satisfies
\begin{align*}
{\phi}'(s) > 0,
\end{align*}
and 
\begin{align*}
{\phi}''(s) - {\phi}'(s)^2 - {\phi}'(s) > 0.
\end{align*}
\end{lemma}
\begin{proof}
 This follows easily from Proposition \ref{choosephi}.
\end{proof}
With $\phi(s)$ chosen as above, we let 
$$\epsilon_1 = \mbox{min}\{ \phi'(s) \},$$ and
$$\epsilon_2 = \mbox{min} \{ \phi''(s) - \phi'(s)^2 + \phi'(s) \}.$$  
From the inequality (\ref{neggoodone}), we have
\begin{align*}
C \geq 
%\frac{1}{2v} T_{k-1}^{ij} \Big( \partial_i \partial_j
%g^{lm} + 2 \partial_l \Gamma^m_{ij} \Big) u_l u_m   \\
 \epsilon_2  T_{k-1}^{ij} u_i u_j
 + T_{k-1}^{ij} \left( 
 R_{iljm} \frac{u_l u_m}{v} + \epsilon_1 \frac{| \nabla u|^2}{2} 
\delta_{ij} + \phi'(u) S_{ij}
+ \frac{u_l}{v} \partial_l S_{ij} \right).  
\end{align*}
The proof then procedes exactly as before. 
\end{proof}
   We note that our method above for obtaining 
the $C^2$ estimate fails for equation (\ref{specific}),
since the dominating term in the inequality
(\ref{coolineq}) now has the wrong sign. 
%%%%%%%%%%%%%%%%%%%%%%%%%%%%%%%%%%%%%%%%%%%%%%%%%%%%%%%%%%%%%%%%%
\section{Monge-Amp\`ere equation in conformal geometry}
\label{conformalsection}
 In this section we restrict our attention to 
$k=n$, the determinant, and we consider more 
generally:
\begin{align}
\label{deteqn}
{\mbox{det}}^{1/n} \left( \nabla^2 u + du \otimes du - \frac{|\nabla u|^2}{2}g
+ S \right) =  e^{-2u},
\end{align}
where $S \in \Gamma_n^+$ is a positive definite symmetric tensor.

\subsection{Proof of Theorem \ref{detthm}}
 We begin by proving a Harnack 
inequality for solutions of (\ref{deteqn}).
\begin{proposition}
\label{aaaa}
Let $u$ be a $C^2$ solution of (\ref{deteqn}). 
If $\lambda_{max}(S) D^2 < \frac{\pi^2}{2}$, then 
\begin{align}
\label{harnack}
 2 \mbox{\em{ log}} \left( 
\mbox{\em{cos}} \left( D\sqrt{ \lambda_{max} (S)/2}
\right) \right)
+ \mbox{\em{ sup }} u <
 \mbox{\em{inf }} u,
\end{align}
where $\lambda_{max}(S)$ denotes the maximum 
eigenvalue of $S$ on $N$ and $D$ is the diameter. 
\end{proposition}
\begin{proof}
 In order to prove this, it 
is convenient to write the equation (\ref{deteqn}) 
in slighty different form. Writing $e^u = v^2$,
with $v>0$, we see that $v$ solves the 
equation
\begin{align}
\label{geodet}
{\mbox{det}}^{1/n} \left( v \nabla^2 v + dv \otimes dv 
- |\nabla v|^2g + \frac{1}{2}v^2 S \right) =  \frac{v^{-2}}{2},
\end{align}
As seen in Section \ref{ellip}, we must have 
\begin{align} 
v \nabla^2 v + dv \otimes dv - |\nabla v|^2g + \frac{1}{2}v^2 S 
 \in \Gamma_n^+,
\end{align}
and therefore since $v>0$, 
\begin{align} 
\label{convexity2}
 \nabla^2 v + \frac{1}{2}vS \in \Gamma_n^+.
\end{align}
Next choose $p \in N$ such that $v(p) = \mbox{sup }v$
and $q \in N$ such that $v(q) = \mbox{inf }v$.
Let $\gamma: [0, d(p,q)] \rightarrow N$ be a unit speed 
minimal geodesic such that $\gamma(0) = p$ and $\gamma(d(p,q)) = q$.  
Letting $\overline{v}$ denote the restriction of $v$ to $\gamma$, we have
\begin{align*}
\overline{v}''(t) + \frac{1}{2} S( \overset{\cdot}{\gamma}(t), 
\overset{\cdot}{\gamma}(t)) \overline{v}(t)> 0,
\end{align*}
therefore 
\begin{align*}
\overline{v}''(t) + \frac{1}{2} \lambda_{max}(S) \overline{v}(t)> 0,
\end{align*}
Let $M = v(p) = \mbox{sup }v$, and 
$\alpha =  \lambda_{max}(S) /2$. 
Then $w(t) = M \mbox{cos}( \sqrt{\alpha} \cdot t)$ 
satisfies 
\begin{align*}
w''(t) + \alpha w(t) = 0, \  w(0) = M = \overline{v}(0), \ w'(0) = 0 = 
\overline{v}'(0).
\end{align*}
If we let $h(t) = (\overline{v}/w)(t)$, then it is easy to verify
that $h$ satisfies the inequality
\begin{align*}
h'' > 2 \sqrt{\alpha} \mbox{ tan}( \sqrt{\alpha}
\cdot t) h', 
\end{align*}
for $\sqrt{\alpha} \cdot t < \pi / 2$.
Integrating this, and using the boundary 
condition $h'(0) = 0$, we find that 
$h'(t) > 0$ for $t>0$. Since $h(0) = 1$, we 
conclude that  
$ \overline{v}(t) > w(t)$ as long as 
$0 <\sqrt{\alpha} \cdot t < \pi / 2$.
Evaluating this at the endpoint $q$, we have 
\begin{align*}
v(q) = \mbox{inf }v > \mbox{sup }v \cdot
\mbox{cos} ( \sqrt{\alpha} \cdot d(p,q))
\geq  \mbox{sup }v \cdot
\mbox{cos} ( \sqrt{\alpha} \cdot D),
\end{align*}
that is,  
\begin{align*}
\mbox{sup }v < 
\left(\mbox{cos}( D\sqrt{ \lambda_{max} (S)/2}) \right)^{-1}
\mbox{inf }v, 
\end{align*}
which implies the stated inequality for $u$. 
\end{proof}
\begin{proposition}
\label{bbbb}
Let $u$ be a solution of (\ref{deteqn}), then there 
exist constants $ \underline{\delta} \leq \overline{\delta}$ 
depending only upon $g, S $ so that 
$\mbox{ \em{sup }} u >  \underline{\delta}$, and
$\mbox{ \em{inf }} u <  \underline{\delta}.$
\end{proposition}
\begin{proof} 
 This follows from the proof of Proposition \ref{easyC0}, but
since we have $e^{-2u}$ instead of $e^{u}$, the 
inequalities are reversed. 
\end{proof}
\noindent
Combining Propositions \ref{aaaa} and \ref{bbbb}, we 
obtain the $C^0$ estimate:
\begin{theorem} 
Let $u$ be a solution of (\ref{deteqn}). If 
\begin{align}
\label{condition}
\lambda_{max}(S)D^2 < \frac{\pi^2}{2},
\end{align}
then there exist a constant $C$ depending only upon $g, S$ so that
$ |u| \leq C$. 
\end{theorem}

  Next, using this a priori estimate, we give a fixed point 
argument to prove the existence of a solution to $\ref{deteqn}$.

\begin{lemma} 
\label{firstex}
If $S \in \Gamma_n^+$ satisfies (\ref{condition}), 
and $0 < f(x) \in C^{\infty}(N)$ then the equation 
\begin{align}
\label{expeqn}
\mbox{ \em{det}}^{1/n} \left( \nabla^2 u + du \otimes du - \frac{|\nabla u|^2}{2}g
+ S \right) =  f(x)e^{- \langle u \rangle}
\end{align}
admits a unique solution $u \in C^{\infty}(N)$
where $ \langle u \rangle = \int_N u \ dvol_g.$ 
\end{lemma}
\begin{proof}
We use the continuity method. For $t \in [0,1]$ we consider 
the equation 
\begin{align}
\label{conteqn}
F_t(u_t) =  \mbox{ {det}} \left( \bar{\nabla}^2 u_t  \right) -
f(x)^{nt} e^{-n \langle u_t \rangle},
\end{align}
where
\begin{align}
\bar{\nabla}^2 u_t \equiv \nabla^2 u_t + du_t \otimes du_t - 
\frac{|\nabla u_t|^2}{2}g + (1-t) \lambda_{max}(S)g + t S. 
\end{align}
Letting $A^{2, \alpha}_t(N) = \{ u \in C^{2, \alpha}(N) : \bar{\nabla}^2 u_t
\in \Gamma_n^+ \}$, we know from Section \ref{ellip} that 
a solution necessarily lies in $A^{2, \alpha}_t(N)$, and 
we claim that the map $F_t : A^{2, \alpha}_t(N) \rightarrow
C^{\alpha}(N)$ is locally invertible at a solution. 
From (\ref{lin}) above we see that the linearized operator is
\begin{align}
 F_t'(u_t)(h) =
T_{n-1} ( \bar{\nabla}^2 u_t )^{ij}
( \nabla^2_{\tilde{g}_t} h)_{ij} 
+n f(x)^{nt}e^{-n\langle u_t \rangle} \langle h \rangle,
\end{align}
where $\tilde{g}_t = e^{-2u_t}g$.
The coefficient matrix $T_{n-1} ( \bar{\nabla}^2 u_t)$
is positive definite, but there is a slight difficulty 
due to the fact that the linearized operator is not 
formally self-adjoint. Nevertheless, it is still 
invertible. This was proved for Monge-Amp\`ere equations 
in \cite{Delanoe2}, and the proof given there is applicable 
in this case. Local invertibility of $F_t$ follows from 
the implicit function theorem (see \cite{GT}).  

  Let $u_t \in C^{2, \alpha}(N)$ be a solution of (\ref{conteqn}).
The matrix $S_t \equiv (1-t) \lambda_{max}(S)g + t S$
satifies the condition (\ref{condition}) for all $t \in [0,1]$,
therefore we have that $u_t$ satisfies the Harnack 
inequality (\ref{harnack}).
Let $q \in N$ be a point where $u_t$ attains a
global minimum. We have 
\begin{align}
\mbox{det}^{1/n}(S_t) \leq f(q) e^{- \langle u_t \rangle},
\end{align} 
which implies $ \langle u_t \rangle < C$. By also considering 
a maximum of $u_t$, we obtain the estimate 
$| \langle u_t \rangle | \leq C$.
Combining this with the Harnack inequality, we obtain 
an a priori $L^{\infty}$ estimate on $u_t$, independent 
of $t$. From the work in Sections \ref{C1} and \ref{C2},
and Evans-Krylov, we obtain an a priori bound on the 
$C^{2, \alpha}$ norm of $u_t$, independent of $t$ for 
some $\alpha \in (0,1)$. Standard
elliptic theory gives a uniform bound on the $C^{k, \alpha}$
norm for each $k \geq 3$. 

We consider the equation $F_0(u_0) = 0$:
\begin{align}
\label{por}
\mbox{det}^{1/n} \left( \nabla^2 u_0 + du_0 \otimes du_0 - 
\frac{|\nabla u_0|^2}{2}g + \lambda_{max}(S)g \right) =
e^{- \langle u_0 \rangle}.
\end{align}
Let $u_0$ be any solution to (\ref{por}). 
As before, by going to a maximum and minimum of 
$u_0$, we find that $e^{- \langle u_0 \rangle} = \lambda_{max}(S)$. 
Then from the arithmetic-geometric inequality, we have 
\begin{align*}
\lambda_{max}(S) = e^{- \langle u_0 \rangle} \leq \frac{1}{n}\Delta u_0 + 
\frac{2-n}{2n}|\nabla u_0|^2 + \lambda_{max}(S).
\end{align*}
We conclude that $\Delta u_0 \geq 0$, which implies 
$u_0 = \mbox{constant}$. The existence 
of a solution at $t=1$ now follows from the continuity 
method. 

  It remains to prove the uniqueness at $t=1$. 
To see this, if we have 2 distinct solutions $u_1$ and 
$v_1$ at $t=1$, we may run the continuity 
method in reverse. From uniqueness at $t=0$, the 
paths we obtain must hit at some time 
$t_0 \in [0,1)$. But since the linearization is invertible
at $t_0$, this contradicts local invertibility.
%To this end, let 
%$u_t$ and $v_t$ both satisfy $F_t(u_t) = F_t(v_t) = 0. $
%Consider the family  $(1-s)u_t + 
%s v_t$ for $s \in [0,1]$.
%From section \ref{ellip} we know that our equation 
%is elliptic along this path so we 
%conclude that 
%\begin{align} 
%L (v-u ) = e^{ - \langle v \rangle} - e^{ - \langle u \rangle} ,
%\end{align}
%where $L$ is a linear second order elliptic operator with no 
%zeroth order terms. Since $L(v-u) = \mbox{constant}$, from 
%Hopf's maximum principle we have $u - v = \mbox{constant}$.
%Returning to the equation, we see that 
%$\langle u \rangle = \langle v \rangle$, and therefore 
%$u = v$. 
\end{proof}
\begin{theorem} 
\label{extnz}
If $A_g \in \Gamma_n^+$ satisfies 
$\lambda_{max}(A_g) D^2 < \frac{\pi^2}{2}$ then the equation 
\begin{align}
\label{gagagoo}
\mbox{ \em{det}}^{1/n} \left( \nabla^2 u + du \otimes du -
  \frac{|\nabla u|^2}{2}g + A_g \right) =  e^{- 2u}
\end{align}
admits a solution $u \in C^{\infty}(N)$. 
\end{theorem}
\begin{proof}
We will employ a fixed point argument 
using the existence and uniqueness of solutions 
to (\ref{expeqn}) in Lemma~\ref{firstex}. 
For $t \in [0,1]$, $ \alpha \in (0,1)$, and 
$u \in C^{2, \alpha}(N)$, let $u_t = H(u,t)$ denote 
the {\em{unique}} solution in $C^{2, \alpha}(N)$
of the equation:
\begin{align}
\label{jacko}
\mbox{ {det}}^{1/n} \left( \nabla^2 u_t + du_t \otimes du_t 
- \frac{|\nabla u_t|^2}{2}g + t A_g + (1-t) \lambda_{max}(A_g)g
\right) = e^{-2tu} e^{- \langle u_t \rangle}
\end{align}
It is easy to show that for each $u \in C^{2, \alpha}(N)$,
the mapping $H(u,t) : [0,1] \rightarrow C^{2, \alpha}(N)$
is uniformly continuous in $t$, and  
we also claim that for each $t \in [0,1]$, $H(u,t) : C^{2, \alpha}(N) 
\rightarrow C^{2, \alpha}(N)$ is a compact operator.
For a bounded subset of $C^{2, \alpha}(N)$, the right hand 
side is bounded in $C^{2, \alpha}(N)$. From the proof 
of Lemma~\ref{firstex}, solutions are bounded in 
$C^{3,\alpha}(N)$. Since $C^{3, \alpha} \subset C^{2, \alpha}$
is a compact embedding, the claim follows. 

We next show that for all $t \in [0,1]$, solutions
of the equation $u = H(u,t)$ satisfy an 
a priori bound $ \Vert u \Vert_{C^{2, \alpha}(N)} < C$. 
As in the proof of Lemma \ref{firstex}, we need only obtain 
an $L^{\infty}$ estimate. 

To this end, let $u_t \in C^{2, \alpha}(N)$ be a fixed
point $H( u_t, t) = u_t$, and $q \in N$ be a point 
where $u_t$ attains a global minimum. Then we have at $q$, 
\begin{align*}
 \mbox{{det}}^{1/n}(A_g(q)) \leq e^{ -\langle u_t \rangle} 
 e^{-2tu_t(q)},
\end{align*}
which implies 
\begin{align}
\label{zap}
 C_1 \leq - \langle u_t \rangle - 2t \ \mbox{inf }u_t,
\end{align}
for some constant $C_1$, and we obtain the estimate
\begin{align*}
(Vol(N) + 2t)\mbox{inf }u_t  \leq C_1.
\end{align*}
Similary by considering a maximum of $u_t$ we obtain
\begin{align*}
(Vol(N) + 2t)\mbox{sup }u_t \geq C_2,
\end{align*}
for some constant $C_2$.
These estimates, coupled with the Harnack inequality 
in Proposition~\ref{harnack}, imply the desired uniform $L^{\infty}$
estimate.

 As already seen in the proof of Lemma \ref{firstex},
we have that $H(u, 0) \equiv C$ for all $u \in C^{2, \alpha}$,
where $C$ is some constant. We may then apply a fixed point theorem 
of Berger \cite[Theorem 5.4.14, page 270]{Berger}:
\begin{proposition}
Let $H(x,t)$ be a one-parameter family of compact 
operators defined on a Banach space $X$ for $t \in [0,1]$, with 
$H(x,t)$ uniformly continuous in $t$ for fixed $x \in X$. 
Furthermore, suppose that every solution of $x = H(x,t)$
for some $t \in [0,1]$, is contained in the fixed 
open ball $ \sigma = \{ x | \Vert x \Vert < M \}$.
Then, assuming $H(x,0) \equiv 0$, the compact 
operator $H(x,1)$ has a fixed point $x \in \Sigma$. 
\end{proposition}
Letting $X = C^{2, \alpha}(N)$, we find a fixed 
point $u \in C^{2, \alpha}(N)$ at $t=1$.
Standard regularity theory then implies that $u \in C^{\infty}(N)$.
Adding a constant if necessary, we obtain 
a solution to (\ref{gagagoo}).
\end{proof}

 To finish the proof of Theorem \ref{detthm}, 
if $\sigma([g]) < \frac{\pi^2}{2}$, then there exists a metric 
$\bar{g} \in [g]$ with $\lambda_{max}(A_{\bar{g}}) D^2 < \frac{\pi^2}{2}$. 
The existence of a conformal metric $\tilde{g}$ with 
$ \mbox{det}(A_{\tilde{g}}) = 1$ follows 
from Theorem \ref{extnz}, and the compactness of the 
space of such solutions was also demonstrated in the proof of 
Theorem \ref{extnz}. 

\subsection{Examples}
In this section we examine some simple cases, and we refer
the reader to \cite{Petersen} for details.\\
\\
$ \bullet \ \ (S^n, g = \mbox{round metric}): 
Ric = (n-1)g, D = \pi$, and
\begin{align*}
\lambda_{max}(A_{g}) D^2 = 
\pi^2/2.
\end{align*}
If $\sigma(S^n, g) < \frac{\pi^2}{2}$, then Theorem~\ref{detthm}
would imply that the space of solutions of (\ref{zoob}) 
is compact. But compactness cannot hold 
in this case since $S^n$ has a non-compact group of conformal 
transformations, and the orbit of the standard
metric gives rise to a non-compact family 
of solutions of (\ref{zoob}). 
Therefore $\sigma(S^n, g) = \frac{\pi^2}{2}$.
\\
\\
$ \bullet \ \ (RP^n, g = \mbox{standard metric}):
Ric = (n-1)g, D = \pi / 2$, and
\begin{align*} 
\lambda_{max}(A_{g})D^2 = \pi^2/8 < \pi^2/2. 
\end{align*}
From \cite{Jeff2}, we know that the standard metric 
on $RP^n$ is the unique solution in its conformal class 
of ({\ref{zoob}), but this shows that 
conformal classes on $RP^n$ in a large neighborhood of
the standard metric have compactness. 
\\
\\
$ \bullet \ \ (CP^m, g=\mbox{Fubini-Study}): 
Ric = (2m+2)g, D = \pi / 2$, and 
\begin{align*}
\lambda_{max}(A_{g}) D^2 = 
\frac{m+1}{2m-1}\frac{\pi^2}{4} < \pi^2 / 2. 
\end{align*}
In this case, we do not know if the Fubini-Study metric is 
the unique solution in its conformal class to (\ref{zoob})
since it is not locally conformally flat, but the above 
shows that the space of solutions is compact, and also for 
conformal classes on $CP^m$ in a large neighborhood of Fubini-Study. 
\bibliography{references}
\small{\textsc{Department of Mathematics, Massachusetts Institute
of Technology, Cambridge, MA 02139}}\\
{\em{E-mail Address:}} \ {\texttt{jeffv@math.mit.edu}}
\end{document}